\documentclass[11pt,twoside]{amsart}
\newcommand{\tab}{\,\,\,\,}
\newcommand{\be}{\begin{equation}}
\newcommand{\ee}{\end{equation}}
\newcommand{\bea} {\begin{eqnarray}}
\newcommand{\nbea} {\begin{eqnarray*}}
\newcommand{\eea} {\end{eqnarray}}
\newcommand{\del}{{\partial}}
\newcommand{\RR}{\mathbb{R}}
\newcommand{\eps}{\epsilon}
\newcommand{\lam}{\lambda}
\newcommand{\overu}{{u}}
\newcommand{\overv}{{v}}
\newcommand \sign {\text{sgn}}
\newtheorem{definition}{Definition}[section]
\newtheorem{lemma}[definition]{Lemma}
\newtheorem{theorem}[definition]{Theorem}

\numberwithin{equation}{section}

\begin{document}
\title[Validity of Chapman-Enskog Expansions]
{On the Validity of Chapman-Enskog Expansions for Shock Waves with Small Strength}

\author[Bedjaoui, Klingenberg, and L{\tiny e}Floch]
{
Nabil Bedjaoui$^{1,3}$, Christian Klingenberg$^2$,
\\
and
\\
Philippe G. L{\tiny e}Floch$^3$}
\thanks{$^1$ INSSET, Universit\'e de Picardie, 48 rue Raspail, 02109
Saint-Quentin, France.  E-mail: {\sl bedjaoui@cmap.polytechnique.fr}}
\thanks{$^2$ Applied Mathematics
Department, W\"urzburg University, Am Hubland, 97074 W\"urzburg,
Germany. E-mail: {\sl klingen@mathematik.uni-wuerzburg.de}}
\thanks{$^3$  
 Current affiliation: 
 Laboratoire Jacques-Lions Lions \& Centre National de la Recherche Scientifique (CNRS), 
Ê      Universit\'e Pierre et Marie Curie (Paris 6), 
 Ê     4 Place Jussieu, 75252 Paris, France. 
E-mail: {\sl pgLeFloch@gmail.com.} 
\textit{Key Words:} conservation law, hyperbolic, shock wave,
traveling wave, relaxation, diffusion, Chapman-Enskog expansion.
{Published in:} {\tt Portugal. Math. 61 (2004). 479--499.}
}
\begin{abstract}
We justify a Chapman-Enskog expansion for discontinuous solutions
of hyperbolic conservation laws
containing shock waves with {\sl small} strength. Precisely,
we establish pointwise uniform estimates for the difference
between the traveling waves of a relaxation model
and the traveling waves of the corresponding diffusive equations
determined by a Chapman-Enskog expansion procedure to first- or second-order.
\end{abstract}
\maketitle
\newpage 


\section{Introduction}

We consider scalar conservation laws of the form 
\bea \del_t u +
\del_x f(u) =0, \tab u=u(x,t) \in \RR, \, t >0, 
\label{1.1} \eea 
where the
flux-function $f: \RR \to \RR$ is a given, smooth mapping. It is
well-known that initially smooth solutions of (1.1) develop
singularities in finite time and that weak solutions satisfying
(1.1) in the sense of distributions together with a suitable
entropy condition must be sought. For instance, when the initial
data have bounded variation, the Cauchy problem for (1.1) admits a
unique entropy solution in the class of bounded functions with
bounded variation. (See, for instance, \cite{LeFloch}.) In the
present paper, we are primarily interested in shock waves of
(1.1), i.e. step-functions propagating at constant speed.

Entropy solutions of (1.1) can be obtained as limits of
diffusion or relaxation models. For instance, under the
sub-characteristic condition \cite{Liu}
\begin{equation}
\sup| f'(u)| < a,
\label{1.2}
\end{equation}
and when the relaxation parameter $\eps>0$ tends to zero
it is not difficult to check that solutions of
\bea
\begin{array}{l}
\displaystyle \del_t u_\eps + \del_x v_\eps = 0,
\\
\displaystyle  \del_t v_\eps + a^2 \, \del_x u_\eps = {1 \over \eps} \,
\bigl(f(u_\eps) - v_\eps \bigr), 
\end{array}
\label{1.3} 
\eea 
converge toward entropy solutions of (1.1). More
precisely, the first component $u:=\lim_{\eps \to 0} u^\eps$ is an
entropy solution of (1.1) and $f(u):= \lim_{\eps \to 0} v^\eps$ is
the corresponding flux. See, for instance, Natalini \cite{Natalini}
and the references therein for a review and references.

The Chapman-Enskog approach \cite{ChapmanCowling} allows one to
approximate (to ``first-order'') the relaxation model (1.3) by a
diffusion equation ((1.4) below). More generally, it provides a
natural connection between the kinetic description of gas dynamics
and the macroscopic description of continuum mechanics. The
Chapman-Enskog expansion and its variants have received a lot of
attention, from many different perspectives. For recent works on
relaxation models like (1.3), Chapman-Enskog expansions, and
related matters we refer to Liu \cite{Liu}, Caflisch and Liu
\cite{CaflischLiu}, Szepessy \cite{Szepessy}, Natalini
\cite{Natalini}, Mascia and Natalini \cite{MN}, Slemrod
\cite{Slemrod}, Jin and Slemrod \cite{JinSlemrod}, Klingenberg and
al. \cite{KlingenbergLuZhao}, and the many references therein.

Our goal in this paper is to initiate the investigation
of the validity of the Chapman-Enskog expansion for discontinuous solutions
containing shock waves.  This expansion is described in the literature
for solutions which are sufficiently smooth, and
it is not a priori clear that such a formal procedure could still be
valid for {\sl discontinuous} solutions.
This issue does not seem to have received the attention it deserves, however.
Note first that, by the second equation in (1.3), we formally have
\be
\begin{split}
v_\eps
& = f(u_\eps) - \eps \, \bigl(\del_t v_\eps + a^2 \, \del_x u_\eps\bigr)
\nonumber
\\
&  =  f(u_\eps) - \eps \, \bigl(\del_t f(u_\eps) + a^2 \, \del_x u
_\eps\bigr) + O(\eps^2) \\
&  =  f(u_\eps) - \eps \, \bigl(-f'(u_\eps) \, \del_x f(u_\eps) +
a^2 \, \del_x u _\eps\bigr) + O(\eps^2), \nonumber
\end{split}
\ee
as long as second-order derivatives of the solution remain uniformly bounded in $\eps$.
Keeping first-order terms only, we arrive at the diffusion equation
\be
\del_t u_\eps + \del_x f(u_\eps) =  \eps \, \del_x \bigl( (a^2 -
f'(u_\eps)^2 ) \,
\del_x u _\eps\bigr).
\label{1.4}
\ee
This expansion can be continued at higher-order to provide, for smooth
solutions of (1.3), an approximation with higher accuracy.
When solutions of (1.3) cease to be smooth and the gradient $\del_x
u_\eps$ becomes
large, the terms collected in $O(\eps^2)$ above are clearly no longer
negligible
in a neighborhood of jumps. The validity of the first-order approximation (1.4), as well
as higher-order expansions in powers of $\epsilon$, becomes
questionable.

The present paper is motivated by earlier results
by Goodman and Majda \cite{GoodmanMajda} (validity of the equivalent
equation associated with a difference scheme), Hou and LeFloch
\cite{HouLeFloch}
(difference schemes in nonconservative form),
and Hayes and LeFloch \cite{HayesLeFloch} (diffusive-dispersive schemes
to compute nonclassical entropy solutions).
In these three papers, the validity of
an asymptotic method is investigated for {\sl discontinuous} solutions, by
restricting attention to {\sl shock waves with sufficiently
small strength.}
This is the point of view we will adopt and, in the present paper,
we provide a rigorous justification of the validity of the
Chapman-Enskog expansion for solutions containing shocks with small strength.

Specifically, restricting attention to traveling wave solutions of the
relaxation model (1.3), the first-order approximation (1.4),
and the associated second-order approximation (see Section 2 below),
we establish several pointwise, uniform estimates which show that
the first- and the second-order approximations approach closely the
shock wave solutions of (1.3) with sufficiently small strength.
See Theorem 3.2 (for Burgers equation), Theorem 4.2 (general conservation laws),
and
Theorem 5.1 (generalization to second-order approximation).
In the last section of the paper, we discuss whether our results are expected to
generalize to higher-order approximations.


\section{Formal Chapman-Enskog expansions}

\subsection{Expanding $v_\eps$ only} In this section we will discuss two variants
to derive a formal Chapman-Enskog expansion for (1.3), at any
order. We begin by plugging the expansion $ v = \sum_{k=0}^\infty
\eps^k \, v_k$ into (1.3) while keeping $u$ fixed. We obtain \bea
\begin{array}{l}
\displaystyle \del_t u + \sum_{k=0}^\infty \eps^k \,  \del_x v_k = 0,
\\
\displaystyle  \sum_{k=0}^\infty \eps^k \, \del_t  v_k
   + a^2 \, \del_x u = {f(u) \over \eps} - \sum_{k=0}^\infty \eps^{k-1} \,
v_k.
\end{array}
\nonumber \eea The second identity above yields
 \bea
\begin{array}{l} \nonumber
f(u) = v_0,
\\
\del_t v_0 + a^2 \, \del_x u = - v_1,
\\
\del_t v_k = - v_{k+1}, \quad k\geq 1,
\end{array}
\label{1}
\eea
which determines $v_0 = f(u)$ and, for $k \geq 1$, $v_k = {(-1)}^k
\del_t^{k-1} \big( \del_t f(u) + a^2 \, \del_x u\big)$, while the
function $u$ is found to satisfy
\be
\del_t u + \del_x f(u) = - \del_x \sum_{k=1}^\infty {(-\eps)}^k \,
\del_t^{k-1} \big( \del_t f(u) + a^2 \, \del_x u\big).
\label{2}
\ee

For instance, to first order we find
\be
\del_t u + \del_x f(u) = \eps \del_x
  \big( \del_t f(u) + a^2 \, \del_x u\big),
\label{3}
\ee
and to second order
\be
\del_t u + \del_x f(u) = \eps \del_x
  \big( \del_t f(u) + a^2 \, \del_x u\big)
-\eps^2\del_{xt}\big( \del_t f(u) + a^2 \, \del_x u\big).
\label{4}
\ee

The corresponding traveling wave equation satisfied by solutions
of the form
$$
u(x, t ) = u(\xi), \qquad  \xi := (x -\lam \, t)/\eps
$$
read
\be
- \lam \, u' + f(u)' = \sum_{k=1}^\infty \lam^{k-1} \big( (-\lam \,
f'(u) + a^2)
\, u' \big)^{(k)}.
\label{M1.3}
\ee
To first order the traveling wave equation is
\be
-\lam \, u' + f(u)' = \big( (-\lam \, f'(u) + a^2) \, u' \big)'
\ee
and to second order
  \be
-\lam \, u' + f(u)' = \big( (-\lam \, f'(u) + a^2) \, u' \big)'
            + \lam \big( (-\lam \, f'(u) + a^2) \, u' \big)''.
\label{M1.4}
\ee


\subsection{Expanding both $u_\eps$ and $v_\eps$} One can also expand
both $u^\eps$ and $v^\eps$, as follows:
\bea
\begin{array}{l}
\displaystyle u_\eps = u_0 + \eps \overu_1 + ... = u_0 + \sum_{k=1}^{\infty}
\eps^k\, \overu_k,
\nonumber
\\
\\
\displaystyle v_\eps = v_0 + \eps v_1 + ... =
v_0 + \sum_{k=1}^{\infty}
\eps^k\, \overv_k.
\nonumber
\end{array}
\eea
The solution at $k^{\mbox{th}}$-order is defined by
\be
\widetilde u_k
: = u_0 + \eps \, \overu_1 +... + \eps^k \, \overu_k.
\ee
We also set
\be
\widetilde v_k : = v_0 + \eps \, \overv_1 +... + \eps^k \, \overv_k.
\ee

To first order, one can write (1.3) as \bea
\begin{array}{l}
\displaystyle \nonumber
  \del_t u_0 + \eps \,
\del_t u_1 + \del_x v_0 + \eps \, \del_x \overv_1 + O(\eps^2) = 0,
\\
\displaystyle  \nonumber \del_t v_0 + \eps \, \del_t \overv_1 +
a^2 \, ( \del_x u_0 + \eps \, \del_x \overu_1) + O(\eps^2)= { 1 \over \eps}
\big( f(u_0) + \eps \, f'(u_0) \, \overu_1 - v_0 - \eps \,
\overv_1 \big) + O(\eps),
\end{array}
\eea
which yields the following equations:
\bea
\begin{array}{l} \nonumber
\displaystyle f(u_0) - v_0 = 0
\\
  \del_t u_0 + \del_x v_0 = 0,
\\
\nonumber\displaystyle \del_t \overu_1 + \del_x \overv_1 = 0,
\\
\nonumber\displaystyle \del_t v_0 + a^2 \, \del_x u_0 =  f'(u_0)
\, \overu_1 - \overv_1.
\end{array}
\eea
Thus \bea
\begin{array}{l}
\displaystyle   \nonumber \del_t u_0 +  \del_x f(u_0) = 0,
   \\
\displaystyle  \nonumber \del_t \overu_1 + \del_x \big( f'(u_0) \,
\overu_1 - \del_t v_0 - a^2 \, \del_x u_0 \big) = 0.
\end{array}
\eea
Therefore, the first-order, Chapman-Enskog expansion leads us to
\be
\begin{split}
\del_t \widetilde u_1 + \del_x f(\widetilde u_1)
& = \del_t (u_0 + \eps \, \overu_1) + \del_x \big(f(u_0) +
\eps \, f'(u_0) \, \overu_1\big)
\\
& = \eps \, (\del_{xt}v_0 + a^2 \del_{xx} u_0)
\\
& = \eps \, \big(a^2 \del_{xx}u_0 - \del_{tt}u_0\big).
\end{split}
\ee
Using that $\del_t u_0 = - \del_x f(u_0)$
we get
\bea
\nonumber
\del_t  \widetilde u_1 + \del_x f(\widetilde u_1)
=
\eps \, \del_x
\big( ( -f'(u_0)^2 + a^2) \, \del_x u_0\big)
+ O(\eps^2).
\eea
Neglecting the terms in $O(\eps^2)$ we may consider that
$\widetilde u_1 = u_0 + \eps \, \overu_1$
is a solution of
\bea
\del_t  \widetilde u_1 + \del_x f(\widetilde u_1) = \eps \,
\del_x \big( ( -f'(\widetilde u_1)^2 + a^2) \, \del_x \widetilde
u_1\big).
\eea

By a similar, but more tedious calculation we can also derive the
diffusive equation at second-order. Using (2.7) and (2.9), we have
\bea \nonumber \del_t \widetilde u_2 + \del_x f(\widetilde  u_2)
&& = \del_t \widetilde u_1 + \eps^2 \del_t \overu_2 + \del_x
\big(f(\widetilde u_1) + \eps^2 \,
f'(\widetilde u_1) \, \overu_2\big)\\
\nonumber
  &&
= \eps\big(a^2 \del_{xx}u_0 - \del_{tt}u_0\big) \,  +
\eps^2\big(\del_t\overu_2 + \del_x \big( f'(\widetilde u_1) \,
\overu_2\big)\big). \eea But, the second order expansion in (1.3)
gives
\bea
\begin{array}{l}
\nonumber\displaystyle \del_t \overu_2 + \del_x \overv_2 = 0,
\\
\nonumber\displaystyle \del_t \overv_1 + a^2 \, \del_x \overu_1 =
f'(\widetilde u_1) \, \overu_2 - \overv_2,
\end{array}
\eea
and we get
\bea
\nonumber
  \del_t \widetilde u_2 + \del_x f(\widetilde u_2)&&
   = \eps\big(a^2 \del_{xx}u_0 - \del_{tt}u_0\big)
+\eps^2\del_x\big( f'(\widetilde u_1) \, \overu_2 -
\overv_2\big)\big)
\\
\nonumber && =\eps\big(a^2 \del_{xx}u_0 - \del_{tt}u_0\big)
+\eps^2\del_x\big( \del_t \overv_1 + a^2 \, \del_x \overu_1\big)
\\
\nonumber && =\eps\big(a^2 \del_{xx}u_0 -
\del_{tt}u_0\big)+\eps^2\big( -\del_{tt} \overu_1 + a^2 \,
\del_{xx} \overu_1\big)
\eea
Finally, since $\widetilde u_1 = u_0 + \eps \overu_1$ we conclude that,
to second order,
\bea
\del_t \widetilde u_2 + \del_x f(\widetilde u_2)  =\eps\big(a^2
\del_{xx} \widetilde u_1 -
\del_{tt} \widetilde u_1\big).
\eea

\

In exactly the same manner we have, for $n\geq 1$, \bea \del_t
\widetilde u_n + \del_x f(\widetilde u_n)  =\eps\big(a^2 \del_{xx}
\widetilde u_{n-1} - \del_{tt} \widetilde u_{n-1}\big), \nonumber
\eea so that \be \del_t \widetilde  u_n + \del_x f( \widetilde
u_n) = \eps \, \big(a^2 \del_{xx} \widetilde u_{n} - \del_{tt}
\widetilde u_{n}\big) + O(\eps^{n+1}). \ee In general, the
$n^{\mbox{th}}$-order equation is obtained by replacing $\del_{tt}
\widetilde u_{n - 1}$ by derivatives with respect to $x$ to obtain
an equation of the form \bea \del_t \widetilde u_n + \del_x
f(\widetilde  u_n) = \sum_{k = 1}^n \eps^k H_k(\widetilde u_n,
\del_x \widetilde u_n,..., \del^{k+1}_x \widetilde u_n). \eea We
will refer to this expansion as  \emph{the Chapman-Enskog
expansion to $n^{\mbox{th}}$ order}.

So let us for instance derive in this fashion the second order
equation satisfied by $\widetilde u_2$. We have first
\bea
\begin{array}{ll}
\displaystyle
\del_{tt} \widetilde u_1 &= \del_t(\del_t \widetilde u_1)
= \del_t \left( -
\del_x f(\widetilde u_1) +
\eps\del_x\big((a^2 - f'(\widetilde u_1)^2)\del_{x} \widetilde
u_1\big)\right) \\
\displaystyle & = - \del_x ( f'(\widetilde u_1)\del_t \widetilde u_1) +
\eps\del_{xt}((a^2 - f'(\widetilde u_1))\del_{x} \widetilde u_1)
\end{array}
\eea
Then setting
$$g_1'(u) = a^2 - f'(u)^2\quad\text{ and }\quad g_2'(u) = (a^2 -
f'(u)^2) f'(u) = g_1'(u) f'(u),
$$
we get
\bea
\begin{array}{ll}
\displaystyle \del_{tt} \widetilde u_1 &=  - \del_x \left(
f'(\widetilde u_1) ( -  f'(\widetilde u_1)\del_x \widetilde u_1 +
\eps \del_{xx} g_1(\widetilde u_1))\right) +\eps
\del_{xxt}g_1(\widetilde u_1)  + O(\eps^2)
\\
  \displaystyle &=
  \del_x ( f'(\widetilde u_1)^2 \del_x \widetilde u_1) - \eps \del_x
(f'(\widetilde u_1) \, \del_{xx} g_1(\widetilde u_1)) +
\eps\del_{xx}(g_1'( \widetilde u_1)\del_t \widetilde u_1)+
O(\eps^2)
\\
\displaystyle &= \del_x ( f'( \widetilde u_1)^2 \del_x  \widetilde
u_1) - \eps \del_x (f'( \widetilde u_1) \, \del_{xx}
g_1(\widetilde u_1)) + \eps \del_{xx}(g_1'( \widetilde u_1)( - f'(
\widetilde u_1)\del_x \widetilde u_1) + O(\eps^2)
\\
\displaystyle &= \del_x ( f'( \widetilde u_1)^2 \del_x  \widetilde u_1)
- \eps \del_x (f'( \widetilde u_1)\del_{xx} g_1( \widetilde u_1))
- \eps
\del_{xxx}g_2( \widetilde u_1) + O(\eps^2).
\end{array}
\nonumber \eea Finally, since $ \widetilde u_1 =  \widetilde u_2 +
O(\eps^2)$, from (2.11) we obtain \bea \del_t  \widetilde u_2 +
\del_x f( \widetilde u_2) = \eps \, \del_{xx} g_1( \widetilde u_2)
+ \eps^2 \, \del_x \left(f'( \widetilde u_2)\del_{xx} g_1(
\widetilde u_2))+ \del_{xx}g_2( \widetilde u_2)\right). \eea
Setting $u = \widetilde u_2$, we can rewrite the last equation in
the form \be
\begin{split}
u_t + f(u)_x = & \eps ((a^2 - f'(u)^2)\, u_x)_x
\\
& + \eps^2 \, \bigl(f'(u) ((a^2 - f'(u)^2)\, u_x)_x\bigr)_x
\\
&+ \eps^2 \, ((a^2 - f'(u)^2) \, f'(u)\,u_x )_{xx}.
\end{split}
\label{M2.1}
\ee
For later reference we record here the traveling wave equation
associated with \eqref{M2.1}
\be
\begin{split}
- \lam \, u' + f(u)'
= & \big( (a^2 - f'(u)^2) \, u'\big)'
\\
& + (f'(u) \, ( (a^2 - f'(u)^2) \, u')')'
\\
& + ((a^2 - f'(u)^2) \, f'(u)\,u' )''.
\end{split}
\label{M2.2}
\ee

\

We arrive at the main issue in this paper : Does the solution
$\widetilde u_n $ of (2.13) converge to some limit $u$ when $n\to
\infty$ and, if so, does this limit satisfy the equation
$$
\del_t u + \del_x f(u)  =\eps\big(a^2 \del_{xx}u - \del_{tt}u\big).
$$
In other word, is this limit $u$ a solution of the relaxation
model (1.3) ? To make such a claim rigorous one would need to
specify in which topology the limit is taken.  As we are
interested in the regime where shocks are present the convergence
in the sense of distributions should be used. We will not address
this problem at this level of general solutions, but will
investigate the important situation of traveling wave solutions,
at least as far as first- and second-order approximations are
concerned.


\section{Burgers equation : Validity of the first-order equations}

We begin, in this section, with the simplest flux function
$f(u) = u^2/2$.
Modulo some rescaling $x\to x-\lam t/\eps$, the traveling wave
solutions $u=u(x)$,
$v=v(x)$ of (1.3) are given by
\bea
\begin{array}{l}
\displaystyle
- \lam \, u' + v' = 0, \\
-\lam \, v' + a^2 \, u' = {u^2 \over 2} - v,
\end{array}
\label{2.1}
\eea
where $\lam$ represents the wave speed. Searching for solutions connecting
left-hand states $u_-$ and $v_-:=f(u_-)$
to right-hand states $u_+$ and $v_+:=f(u_+)$ (so both at
equilibrium), we see that
$$
\lam \, (u_+ - u_-) = v_+ - v_-,
$$
so that the component $u$ is a solution of the single first-order equation
$$
(a^2 - \lam^2) \, u' = {1 \over 2} (u - u_-) (u-u_+).
\label{2.2}
$$
The shock speed is also given by $\lam= (u_+ + u_-)/2$.
Finally,
an easy calculation based on (3.1) yields the following explicit
formula for the solution, say $u=u_*(x)$ of (3.1)
connecting $u_-$ to $u_+$. It exists if and only if  $u_- > u_+$
and then
\be
u_*(x) := u_- - {(u_- - u_+) \over 1 + \exp\left(- {u_- - u_+ \over
2(a^2 - \lam^2)} \, x  \right)}.
\label{2.3}
\ee

It will be useful to introduce the following one-parameter family
of functions \bea \varphi_{ \mu}(x) := u_- - {(u_- - u_+) \over 1
+ \exp\left(- {x (u_- - u_+) \over 2(a^2 - \mu)} \right)}, \qquad
\mu \in \RR \setminus \{a^2\},
\label{2.5}
\eea
in which $\mu$ is a
parameter, not necessarily related to the speed $\lam$.
Clearly, we have
$$
u_* = \varphi_{\lambda^2}.
$$
Note that we have for all $\mu < a^2$, and $x\in \RR$,
$$
  u_+ < \varphi_{ \mu}(x) < u_-
$$

The following estimate in terms of the strength $\delta := (u_- - u_+) $
is
easily derived from (3.3):

\begin{lemma} Given $a>0$ and $0 < h < a^2$ there exist constants
$c$, $C>0$ such that for all $ \mu_1, \mu_2 \in (- a^2+h, a^2 -
h)$ and for all $x \in \RR$
we have
\bea
  |\varphi_{ \mu_1}(x) - \varphi_{\mu_2}(x)| \leq C \, \delta^2 \, |x| \,
|\mu_1 - \mu_2| \, e^{- c\, |x| \, \delta}.
\label{2.6}
\eea
\end{lemma}

\begin{proof}  We  can write
\bea
\begin{array}{l}
\displaystyle   |\varphi_{ \mu_1}(x) - \varphi_{ \mu_2}(x)| =
\delta \, \left|
    {1 \over 1 + \exp\left(- {x (u_- - u_+) \over 2(a^2 - {\mu_2})}
    \right) }
    -
    {1 \over 1 + \exp\left(- {x (u_- - u_+) \over 2(a^2 - { \mu_1})}
    \right) }
    \right| \\
   \displaystyle   \leq {|x| \over 2} \, \delta^2 \left| {1 \over a^2 -
{ \mu_2}} - {1 \over a^2 - {\mu_1}}    \right|
    \, \sup_{x, k}
    {\exp\left(- {x \, \delta \over 2(a^2 - k)}\right)  \over
\left(1 + \exp\left(-
    {x \, \delta \over 2(a^2 - k)} \right) \right)^2}.
\end{array}
\label{2.11} \eea Here, the super bound is taken for $|k|< a^2 -
h$ and $x\in \RR$.

 Then observe that for
$y>0$ we have \bea {\exp\left(- {y\over 2(a^2 - k)} \right) \over
\left(1 + \exp\left(-
    {y\over 2(a^2 - k)} \right)\right)^2 }
    & \leq
     \, \exp\left(- {y\over 2(a^2 - k)}\right)
\nonumber\\
& \leq
     \, \exp\left(- {y\over 2(a^2+ (a^2 - h)) }\right),
\nonumber \eea while for $y<0$ we have

\bea
{\exp\left(- {y\over 2(a^2 - k)} \right) \over \left(1 +
\exp\left(-
    {y\over 2(a^2 - k)} \right) \right)^2}
    & \leq
    {1 \over 1+ \exp\left(- {y\over 2(a^2 - k)}\right)  }
   \nonumber \\
& \leq
     \exp\left( {y\over 2(a^2 - k)}\right)
  \nonumber \\
& \leq \exp\left( {y\over 2(a^2 + (a^2-h))}\right).
\nonumber
\eea
This establishes the desired estimate.
\end{proof}

We are now in position to study the traveling waves of the
first-order equations
obtained by either the approaches in Subsections 2.1 and 2.2:
$$
- \lam \, u' + \left({u^2 \over 2}\right)' =  \bigl( (a^2 - \lam\ u )
\, u'\bigr)'
$$
and
$$
- \lam \, u' + \left({u^2 \over 2}\right)' =  \bigl( (a^2 - u^2 )
\, u'\bigr)',
$$
respectively. Note that they only differ by the diffusion
coefficients in the right-hand sides. After integration, calling
$V_1$ and $W_1$ the corresponding traveling wave solutions, we get
\be (a^2 - \lam\, V_1) \, V_1' = {1 \over 2} \, (V_1 - u_-) \,
(V_1 - u_+) \label{3.6} \ee and \be (a^2 - W_1^2) \, W_1' = {1
\over 2} \, (W_1 - u_-) \, (W_1 - u_+), \label{3.7} \ee
respectively. For uniqueness, since the traveling waves are
invariant by translation, we assume in addition that for example
\bea
u_*(0) = V_1(0) = W_1(0) = {u_- + u_+\over 2}.
\eea

To compare the first-order diffusive traveling waves $W_1$ and $V_1$
with the relaxation
traveling wave $u_*$,
we rely on monotonicity arguments. It is clear that the traveling
waves are monotone, with $V_1'\, , W_1' <0$ and $u_- > V_1(x)\, ,
W_1(x)  > u_+$, so that
setting
$$
\Gamma_- = \min_{[u_+, u_-]} u^2  - b\, \delta, \qquad \Gamma_+ =
\max_{[u_+, u_-]} u^2 + b\, \delta,
$$
where $b>0$ is a sufficiently small constant such that $ \Gamma_+
< a^2$, we find \be
\begin{split}
& (a^2 - \Gamma_-) \, W_1' < {1 \over 2} \, (W_1 - u_-) \, (W_1
- u_+),
\\
& (a^2 - \Gamma_+) \, W_1' > {1 \over 2} \, (W_1 - u_-) \, (W_1 - u_+).
 \label{2.9}
\end{split}
\ee
Therefore, setting
$$
\tilde u = W_1 - \varphi_{\Gamma_-},
$$
after some calculation we find 
\be
2\,(a^2 - \Gamma_-) \, \tilde
u' - \tilde u^2 + \tilde u\, \delta \, {  1 - \exp\left(- {x
\delta\over 2(a^2 - \Gamma_-)} \right)  \over
    1 + \exp\left(- {x \delta\over 2(a^2 - \Gamma_-)} \right) }  < 0.
\label{2.10}
\ee
We have $\tilde u(\pm \infty) =0$. As $x \to \pm
\infty$ the last coefficient in \eqref{2.10} approaches $\pm 1$ and the
function $\tilde u$ satisfies
$$
c \, \tilde u' \pm \tilde u\, \delta + \mbox{H.O.T.} < 0.
$$
So, $\tilde u$ decreases exponentially at infinity while keeping a
constant sign, and we deduce that $\tilde u(x) \neq 0$ for $|x|
\geq M$, for some sufficiently large $M$.

Now, if $\tilde u$ vanishes at some point $x_0$ then, thanks to
the inequality \eqref{2.10}, we deduce that $\tilde u'(x_0) <0$. This
implies that there is at most one point, and thus exactly one
point where $\tilde u$ vanishes, which is by (3.8) $x_0=0$.
Therefore, we have $\sign(x) \, \tilde u(x) < 0$.

A similar analysis
applies to the function $W_1 - \varphi_{\Gamma_+}$  and we obtain
\be \sign(x) \, \varphi_{\Gamma_+}(x) < \sign(x)\, W_1(x) < \sign(x)
\,\varphi_{\Gamma_-}(x), \qquad x \in \RR.
 \label{2.14}
\ee

Concerning the function $V_1$, by defining
$$
\lam_- := \min_{[u_+, u_-]} u, \quad \lam_+ := \max_{[u_+, u_-]}
u,
$$
and
$$
\Lambda_- := \min(\lam \,\lam_-, \lam\,\lam_+)\, - b\, \delta,
\quad \Lambda_+ := \max(\lam\,\lam_-, \lam \,\lam_+)\, + b\,
\delta,
$$
where $b>0$ is a sufficiently small constant such that $\Lambda_+
< a^2$, we obtain in the same manner as above \be \sign(x) \,
\varphi_{\Lambda_+}(x) < \sign(x)\, V_1(x) <\sign(x)
\,\varphi_{\Lambda_-}(x), \qquad x \in \RR.
\label{2.1444}
\ee Note that, for the same reasons, the function $u = u_*$
satisfies also (3.11) and (3.12). Finally, since $|\Gamma_+ -
\Gamma_- |, |\Lambda_+ - \Lambda_-|\ \leq C\, \delta$, we can
combine (3.11) and (3.12) with Lemma 3.1 and conclude:

\begin{theorem} Given two reals $a > M >0$, there are
  constants $c,\,C>0$ so that the following
property holds for all $u_-,u_+ \in [-M,M]$. The uniform distance
between the traveling wave of the relaxation model and the ones of
the first-order diffusive equations derived in Section 2 is of
cubic order, in the sense that \be |V_1(x) - u_*(x)|, \, |W_1(x) -
u_*(x)| \leq C \, \delta^3 \, |x| \, e^{- c\, \delta |x|}, \qquad
x \in \RR. \label{3.11} \ee
\end{theorem}

Note that the estimate is cubic on any compact set but is solely
quadratic in the uniform norm on the real line:
\be
\|V_1 - u_*\|_{L^\infty(\RR)}, \, \|W_1 - u_*\|_{L^\infty(\RR)} \leq
C' \, \delta^2.
\ee


\section{Validity of the first-order expansions}

We extend the result in Section 3 to general, strictly convex flux-functions.
It is well-known that a traveling wave connecting $u_-$ to $u_+$
must satisfy the condition $u_- > u_+$ which we assume from now on.

Set \bea P(u) = f(u) - f(u_-) - \lam \, (u - u_-), \eea and denote
by $u_*$ the solution of the relaxation equation and by $V_1$ and
$W_1$  the first-order traveling waves corresponding to equation
(2.2)(i.e, (2.5)) and to (2.10) respectively. We have \be
\begin{split}
& (a^2 - \lam^2) \, u_*' = P(u_*),
\\
& (a^2 - \lam\, f'(V_1))\, V_1' = P(V_1),
\\
& (a^2 -  f'(W_1)^2)\, W_1' = P(W_1),
\end{split}
\ee
together with the boundary conditions
$$
\lim_{\pm \infty} u_*(x) = \lim_{\pm \infty} V_1(x)
= \lim_{\pm \infty} W_1(x)= u_\pm.
$$

The existence of solutions to these first-order O.D.E.'s can easily be checked,
for instance using the following implicit formula:
$$
F_k(u(x)) - F_k(u(0))  =  x, \quad x \in \RR, \, k = 0, 1, 2,
$$
where
\be
\begin{split}
& F_0'(u) := {(a^2 - \lam^2) \over f(u) - f(u_-) - \lam \, (u -
u_-)}, \tab u \in \RR,
\\
&  F_1'(u) := {(a^2 - \lam\, f'(u)) \over f(u) - f(u_-) - \lam \,
(u - u_-)}, \tab u \in \RR,
\\
& F_2'(u) := {(a^2 - f'(u)^2) \over f(u) - f(u_-) - \lam \, (u -
u_-)}, \tab u \in \RR.
\end{split}
\ee
To ensure uniqueness, we can impose, for example,
\bea
u_*(0) = V_1(0) = W_1(0) = {u_- + u_+\over 2}.
\eea

Now, as was done for Burgers' equation, let us define auxilliary
functions $\varphi_\mu$ as the solutions of
\bea (a^2 - \mu)\,
\varphi_\mu' = P(\varphi_\mu),
\eea
with the same boundary conditions as above. For $ \mu < a^2$ we immediately have
$$
u_+ < \varphi_\mu (x) < u_-, \quad x\in\RR.
$$
Setting  $\delta := (u_- - u_+)$ we get:

\begin{lemma} Suppose that $f$ is a strictly convex flux-function and
$u_- > u_+$. Given $a>0$ and $0 < h < a^2$ there exist constants
$c$, $C>0$ such that, for all $ \mu_1, \mu_2 \in (-a^2+h, a^2-h)$
and for all $x \in \RR$,
\bea
  |\varphi_{\mu_1}(x) - \varphi_{\mu_2}(x)| \leq C \, \delta^2 \, |x| \,
|\mu_1 - \mu_2| \, e^{- c\, |x| \, \delta}.
\eea
\end{lemma}

\begin{proof} Let $\psi$ be the solution of
$$
\psi' = P(\psi) = f(\psi) - f(u_-) - \lambda (\psi - u_-).
$$
We clearly have
$$ \varphi_{\mu}(x) = \psi\left( {x\over a^2 - \mu}\right).$$
Now, we can write
\bea
\nonumber
  |\varphi_{\mu_1}(x) - \varphi_{\mu_2}(x) |&&= \left|\psi\left(
{x\over a^2 - \mu_1}\right) - \psi\left( {x\over a^2 -
\mu_2}\right) \right|\\
\nonumber &&= \left|x\, \psi'(k(x)\, x)\left( {1\over a^2 - \mu_1}
- {1\over a^2 - \mu_2}\right)\right|
\\
&&\leq C\, \left|\mu_1 - \mu_2||x| |P(\psi(k(x)\, x))\right|.
\nonumber \eea Here, $k(x)$ is some real number lying in the
interval  $ \left({1\over a^2 - \mu_1}, {1\over a^2 -
\mu_2}\right)$.

On the other hand we have
$$
|P(\psi(x))|\leq C\, \delta\, |\psi(x) - u_-|\leq C \delta^2.
$$
This implies that \bea |\varphi_{\mu_1}(x) - \varphi_{\mu_2}(x) |
\leq C |\mu_1 - \mu_2|\delta^2 |x|.
 \eea
 The behavior at $\pm
\infty$ is described by
$$
\psi(x) \sim  k_+ e^{ (f'(u_+) - \lam)\, x}, \qquad  x \to +\infty
$$
and
$$
\psi(x) \sim  k_- e^{ (f'(u_-) - \lam)\, x}, \qquad x \to -\infty.
$$
Since the coefficient $k(x)$ is bounded away from 0 and $f'(u_+) -
\lam = c_+ \, \delta$
and $f'(u_-) - \lam = c_- \, \delta$ with $c_+ <0$ and $c_- >0$
(bounded away from zero
since $f$ is strictly convex), this completes the proof.
\end{proof}

Consider now the functions $u_*$, $V_1$ and $W_1$ the solutions of
(4.2). Then, we have:

\begin{theorem} Let $f$ be a strictly convex flux-function, $M>0$ and $a>0$ such that
(1.2) holds in $[-M, M]$. Then there exist  constants $c,\,C>0$ so
that the following inequality holds for all $u_-,u_+ \in [-M,M]$
with $u_-
> u_+$: for all $x \in \RR$ \bea |V_1(x) - u_*(x)| ,\, |W_1(x) -
u_*(x)| \leq C \, \delta^3 \, |x| e^{- c\, \delta |x|}. \eea
\end{theorem}

The proof relies on the following lemma:

\begin{lemma} Suppose that $f$ is a strictly convex flux-function and
$u_- > u_+$.
Assume that $z_+$ and $ z_-$ are the solutions of
$$
z_+' = R_+( z_+),  \qquad z_-' = R_-(z_-), \qquad z_+(0) = z_-(0),
$$
where $R_+=R_+(u)$ and $R_-=R_-(u)$ are any smooth functions
satisfying \bea R_+(u) < R_-(u)<0\qquad \text{for all} \quad
u\in(u_+, u_-).
\eea
Then, the two corresponding curve solutions
cross at $x =0$ only, and \be
\begin{split}
& z_+ > z_- \quad \mbox{ for }  x <0,
\\
& z_+ < z_- \quad \mbox{ for } x >0.
\end{split}
\ee
\end{lemma}

\begin{proof}
If there is $x_0$ such that $ z_+(x_0) = z_-(x_0)$ then thanks to
(4.9),
$$
z_+'(x_0) < z_-'(x_0).
$$
This implies that there cannot be more than one intersection
point. So, $(0, z_+(0))$ is the only interaction point of the two
trajectories, and (4.10) follows as well.
\end{proof}

\

\noindent{\sl Proof of Theorem~4.2.} Setting
$$
\lam_- = \min_{[u_+, u_-]} f'(u),\qquad \lam_+ = \max_{[u_+, u_-]}
f'(u)
$$
and
$$
\Lambda_- = \min(\lambda\,\lambda_-, \lambda\,\lambda_+) - b\,
\delta, \qquad \Lambda_+ = \max(\lambda\,\lambda_-,
\lambda\,\lambda_+) + b\, \delta,
$$
where, $b>0$ is  a sufficiently small constant such that
$\Lambda_+ < a^2$,
 we have
$$
\Lambda_- <  \lam\, f'(u) < \Lambda_+ \quad \text{ and} \quad
\Lambda_- <  \lam^2 < \Lambda_+,
$$
and thus
 \bea
  0< a^2 - \Lambda_+ < a^2 - \lam\, f'(u),\, a^2 - \lam^2
< a^2 - \Lambda_-.
 \eea

Applying Lemma 4.3 we deduce  that \be
\begin{split}
& \varphi_{\Lambda_-}< u_*,\, V_1,\,  < \varphi_{\Lambda_+}\quad x
< 0,
\\
& \varphi_{\Lambda_+}< u_*, \, V_1,\, < \varphi_{\Lambda_-}\quad x
> 0.
\end{split}
\nonumber \ee
 Now, concerning the third equation in (4.2), we set
$$
\Gamma_- = \min_{[u_+, u_-]} f'(u)^2 - b\, \delta\quad \text{
and}\quad \Gamma_+ = \max_{[u_+, u_-]} f'(u)^2 + b\, \delta,
$$
where $b>0$ is sufficiently small such that $\Gamma_+ < a^2$.
 We obtain \bea
  0< a^2 - \Gamma_+ < a^2 - f'(u)^2,\, a^2 - \lam^2
< a^2 - \Gamma_-
 \eea
and, by Lemma~4.3,
 \be
\begin{split}
& \varphi_{\Gamma_-}< u_*,\, W_1,\,  < \varphi_{\Gamma_+}\quad x <
0,
\\
& \varphi_{\Gamma_+}< u_*, \, W_1,\, < \varphi_{\Gamma_-}\quad x
> 0.
\end{split}
\nonumber \ee

 Finally, since $|\Lambda_+ - \Lambda_-|, \, |\Gamma_+ - \Gamma_-| \leq C\,\delta$,
 by applying Lemma~4.1, we obtain (4.8). This
completes the proof of Theorem 4.2. \quad \qed


\section{Validity of a second-order expansion}

Our next objective is to extend the estimate in Theorem 4.2 to
the second-order equation obtained in Subsection 2.1.

We consider the equation (2.6) after integrating it once. The traveling
wave connects $u_-$ to $u_+$, with $u_- > u_+$, and is given by
\bea
P(u) := (-\lam \, f'(u) + a^2) \, u'
            + \lam \, \big( (-\lam \, f'(u) + a^2) \, u' \big)'.
\eea
Defining first- and second- order ODE operators:
$$Q_1 u = (a^2 - \lam f'(u)) u'$$
and
$$
Q_2 u = (a^2 - \lam f'(u)) u' + \lam  \Bigl( (a^2 - \lam f'(u))
u'\Bigr)' = Q_1 u  + \lam (Q_1 u)'.
$$
The solution $u = V_2$ of (2.6) under consideration satisfies
\bea
Q_2 V_2 = P(V_2)
\eea

\begin{theorem} Let $f: \RR \to \RR$ be a strictly convex flux-function
and $M>0$. Then there exist constants $C,c, c_0>0$ so that the following
property holds. For any
 $u_-,u_+ \in [-M,M]$ with $u_- > u_+$ and $0 < \delta = u_- - u_+ < c_0$,
there exists a traveling wave $V_2=V_2(y)$ of (5.2)
connecting $u_-$ to $u_+$. Moreover, this traveling wave
approaches the relaxation traveling wave $u_*$ to fourth-order in
the shock strength, precisely:
\be
|V_2(x) - u_*(x)| \leq
C \, \delta^4 \, |x| \,  e^{ - c\,|x| \, \delta}, \qquad x \in
\RR.
\label{5.4}
\ee
\end{theorem}

The estimate is only cubic in the uniform norm on the whole real line:
\be
\|V_2 - u_*\|_{L^\infty(\RR)} \leq C' \, \delta^3.
\ee

\begin{proof} Setting
$$d_\mu = {\lam \over a^2 - \mu},\quad \text{ and }\quad \gamma_\lam = d_{\lambda^2} =
{\lam \over a^2 - \lam^2},$$
 then  $u_* = \varphi_{\lam^2}$ satisfies
$$
Q_1 u_* = P(u_*) ( 1 + \gamma_\lam (\lam - f'(u_*)) ) = P(u_*) ( 1
- \gamma_\lam P'(u_*) ),
$$
and a simple calculation gives
$$
Q_2 u_* = P(u_*) \bigl( 1 -  \gamma_\lam^2 ( f''(u_*) P(u_*) +
(f'(u_*) - \lam )^2)\bigr) = P(u_*) ( 1 - \gamma_\lam^2
(P\,P')'(u_*) ).
$$
In the same manner, the function $ \varphi_\mu$, that is the solution
of (4.5) satisfies the following equation
$$
Q_1 \varphi_\mu = P(\varphi_\mu) \bigl( 1 + c_\mu  + d_\mu (\lam -
f'(\varphi_\mu)) \bigr)
  = P(\varphi_\mu) \bigl( 1 + c_\mu  - d_\mu P'(\varphi_\mu)\bigr),
$$
where
$$
c_\mu := {\mu - \lam^2\over a^2 - \mu},
$$
and
$$
Q_2 \varphi_\mu = P(\varphi_\mu) \Bigl( 1 +  c_\mu ( 1 + d_\mu
(f'(\varphi_\mu) - \lam)) - d_\mu^2 \bigl( f''(\varphi_\mu)
P(\varphi_\mu) + (f'(\varphi_\mu) - \lam )^2\bigr)\Bigr)
$$
or, equivalently,
$$
Q_2 \varphi_\mu = P(\varphi_\mu) \Bigl( 1 +  c_\mu ( 1 + d_\mu
P'(\varphi_\mu) ) - d_\mu^2 \bigl(  P P'
)'(\varphi_\mu)\bigr)\Bigr).
$$

Now, since $|f'(\varphi_\mu) - \lam| \leq C_0\, \delta$ and
$|f''(\varphi_\mu) P(\varphi_\mu) + (f'(\varphi_\mu) - \lam )^2|
\leq C_0 \delta^2$, then for sufficiently small $\delta$ there
exists a positive constant $C$ such that the following property
holds: by choosing $\mu_+$ and $\mu_-$ in the form
$$
\mu_+ = \lam^2 ( 1 + C \delta^2), \qquad \mu_- = \lam^2 ( 1 -  C
\delta^2),
$$
we obtain
$$
Q_2 \varphi_{\mu_+} =  P(\varphi_{\mu_+})( 1 +
K_+(\varphi_{\mu_+})), \quad \text{where}\quad
K_+(\varphi_{\mu_+})> 0
$$
and
$$
Q_2 \varphi_{\mu_-} =  P(\varphi_{\mu_-})( 1 +
K_-(\varphi_{\mu_-})), \quad \text{where}\quad
K_-(\varphi_{\mu_-})< 0.
$$

Consider the corresponding functions $\varphi_{\mu_+}$ and
$\varphi_{\mu_-}$  and let us use phase
plane argument. The corresponding curves
\be
\begin{split}
& {\mathcal C}_+ : \varphi_{\mu_+}\mapsto (\varphi_{\mu_+},
w_{\mu_+} = Q_1 \varphi_{\mu_+}),
\\
& {\mathcal C}_- : \varphi_{\mu_-}\mapsto (\varphi_{\mu_-},
w_{\mu_-} = Q_1 \varphi_{\mu_-})
\end{split}
\ee
satisfy
\bea
\lam\, l(\varphi_{\mu_+})\,w_{\mu_+}{dw_{\mu_+}\over du} +
w_{\mu_+} = P(\varphi_{\mu_+})( 1 + K_+(\varphi_{\mu_+}))
\eea
and
\bea \lam\, l(\varphi_{\mu_-})\,w_{\mu_-}{dw_{\mu_-}\over du} +
w_{\mu_-} = P(\varphi_{\mu_-})(1 + K_-(\varphi_{\mu_-})),
\eea
where
$$
l(u) := {1\over a^2 - \lam\, f'(u)}.
$$
We claim that the curve $\mathcal C_+$ is ``below'' the curve
${\mathcal C}_-$.

This is true locally near the points $(u_-, 0)$ and $(u_+, 0)$, as
it clear by comparing the tangents to the curves at these points
(using (4.5)). Note that if $\lam = 0$ we have $u = u_*$. We then
distinguish between two cases:

\

{\bf Case 1:} If $\lam >0$, suppose that the two curves issuing
from $(u_-, 0)$, meet for the ``first'' time at some point $(u_0,
w_0)$ with $u_+ < u_0 < u_-$. Then, combining (5.6) and (5.7) at
this point we get
$$
\lam \,l(u_0)\,w_0\,\Bigl({dw_{\mu_+}\over du}(u_0) -
{dw_{\mu_-}\over du}(u_0)\Bigr) = P(u_0)( K_+(u_0) - K_-(u_0)).
$$
This leads to a contradiction, since
$$
w_0 <0, \quad {dw_{\mu_+}\over du}(u_0)\leq {dw_{\mu_-}\over
du}(u_0)\quad \text{ and} \quad P(u_0)( K_+(u_0) - K_-(u_0)) <0.
$$
Consider now the equation (5.2) and let us  study in the phase
plane  the trajectory issuing from $(u_-, 0)$ at $-\infty$.
Comparing the eigenvalues we obtain that the tangent at  this
point lies between those of the reference curves $\mathcal C_+$
and $\mathcal C_-$.
 
In the same manner as before, we obtain that this curve cannot
meet $\mathcal C_+$, nor $\mathcal C_-$, and necessarily converges
to $(u_+, 0)$ as $y \to+\infty$.

\

{\bf Case 2:} If $\lam <0$, we follow the same analysis by
considering the trajectory of (2.5) arriving at $(u_+, 0)$ and the
``last'' intersection point.

\

In both cases, we obtain the existence (and uniqueness) of the
solution of (5.2), denoted by $u = V_2$, and also that its
trajectory called $\mathcal C$ is between $\mathcal C_+$ and
$\mathcal C_-$.

Note that since our equations are autonomous, by choosing $u(0) =
\varphi_{\mu_+}(0) = \varphi_{\mu_-}(0) = (u_- + u_+)/2$, we have
\bea \varphi_{\mu_+} < u <\varphi_{\mu_-}, \qquad x >0 \eea and
\bea \varphi_{\mu_-} < u <\varphi_{\mu_+}\, ,\qquad x <0. \eea
Indeed, from the phase plane analysis, if for some $x_0\in \RR$,
$u(x_0) = \varphi_{\mu_+}(x_0)$ then necessarily $w(x_0)  >
w_{\mu_+}(x_0)$ and then $u'(x_0) > \varphi_{\mu_+}'(x_0)$. This
means that the curves $x\mapsto u(x)= V_2(x)$ and $x\mapsto
\varphi_{\mu_+}(x)$ have only one intersection point, that is $(0,
u(0))$, that satisfies in addition $u'(0) > \varphi_{\mu_+}'(0)$.
We obtain in same manner that the two curves $x\mapsto u(x)$ and
$x\mapsto \varphi_{\mu_-}(x)$ have only one intersection point,
that is $(0, u(0))$, that satisfies in addition $u'(0) <
\varphi_{\mu_-}'(0)$.

Now, using the inequalities (5.8) and (5.9) that are also
satisfied by $u_* = \varphi_{\lam^2}$ ( since $\mu_- < \lam^2 <
\mu_+)$, we can write
$$
\aligned
|u_*(x) - u(x)|
&\leq |\varphi_{\mu_+}(x) - \varphi_{\mu_-}(x)|\\
&\leq |{\mu_+} - {\mu_-}|\delta^2 |x|\, e^{- c |x|\, \delta}\\
&\leq C \delta^4 |x|\, e^{- c |x|\, \delta},
\endaligned
$$
which completes the proof of Theorem 5.1.
\end{proof}


\section{Conclusions}

For the general expansion derived in Subsection 2.2 we now establish an identity
which connects the relaxation equation with its Chapman-Enskog expansion at
{\sl any order} of accuracy. By defining the ODE operator
\bea
Q_n u : = \sum_{k =1 }^n
\lam^{k-1}\left((-\lam\, f'(u)+ a^2)u'\right)^{(k-1)},
\eea
we have:

\begin{theorem}
The traveling wave $u_*$ of the relaxation model satisfies
$$
Q_n u_* = P(u_*)\left( 1 - \gamma_\lam ^n \, R_n (u_*)\right),
$$
where
$\gamma_\lam := \lam / (a^2 - \lam^2)$,
and the remainders $R_n$ are defined by induction :
$$
R_1 : = P', \quad R_{n + 1} : = (P\, R_n)' \quad \text{ for }  n\geq 1.
$$
\end{theorem}

\begin{proof} Note that the ODE operators $Q_n$ satisfy
$$
Q_{n+1}u = Q_1 u + \lam (Q_n u )'.
$$
Now, assume that
$$
Q_n u_* = P(u_*)\left( 1 - \gamma_\lam ^n \, R_n (u_*)\right),
$$
then
$$
Q_{n+1} u_* = P(u_*)\left( 1 - \gamma_\lam  \, P' (u_*)\right) +
\lam \left(P' \left( 1 - \gamma_\lam ^n \, R_n (u_*)\right) -
P(u_*) \gamma_\lam ^n \, R_n' (u_*)\right)u_*'.
$$
But since $u_* ' = {P(u_*)\over a^2 - \lam^2}$ it follows that
$$
Q_{n+1} u_* = P(u_*)\left( 1 - \gamma_\lam^{n+1}  \, (P\, R_n)'(u_*)\right)
=
P(u_*)\left( 1 - \gamma_\lam ^{n+1} \, R_{n+1} (u_*)\right),
$$
which completes the proof.
\end{proof}

Theorem 6.1 provides some indication that, by taking into account
more and more terms in the Chapman-Enskog expansion, the approximating
traveling wave should approach
 the traveling wave equation of the relaxation equation (1.3). For
$n$ large but {\sl fixed} it is conceivable that, denoting $V_n$
the solution of $Q_n u = P(u)$, \be \|V_n - u_*\|_{L^\infty(\RR)}
\leq C_n \, \delta^{n+1}.
 \ee
 However, one
may not be able to let $n \to \infty$ while keeping $\delta$
fixed. In fact, numerical experiments (with Burgers flux)
have revealed that the remainders satisfy only
$$
\|R_n (u_*) \|_{L^\infty} \leq C_n' \, \delta^n,
$$
where the constants $C_n'$ grow exponentially and cannot be
compensated by the factor $\gamma_\lam^n$. One can also easily
check, directly from the definitions, that
$$
\|R_{n} (u_*)\|_{L^\infty} \leq C \delta^n\, n!.
$$

In conclusion, although we successfully established {\sl uniform} error estimates
for first- and second-order models, it is an open problem whether such estimates should
still be valid for higher-order approximations. Theorem 6.1 
indicates that the convergence might hold but, probably, in a {\sl weaker} topology.


\section*{Acknowledgments}
The three authors gratefully acknowledge the support and hospitality of the Isaac
Newton Institute
for Mathematical Sciences at the University of Cambridge, where this research
was performed
during the Semester Program ``Nonlinear Hyperbolic Waves in Phase Dynamics and
Astrophysics'' (Jan.-July 2003), organized by C.M. Dafermos, P.G. LeFloch, and E. Toro. 
PLF was also partially supported by the Centre National de la Recherche Scientifique
(CNRS).
 


\begin{thebibliography}{99}

\bibitem{CaflischLiu} Caflisch R.E. and Liu T.P.,
\textit{Stability of shock waves for the Broadwell equations,}
Comm. Math. Phys. 114 (1988), 103--130.

\bibitem{ChapmanCowling}
Chapman S. and Cowling T.,
{\bf Mathematical theory of non-uniform gases,}
Cambridge University Press, Cambridge, 1970.

\bibitem{GoodmanMajda} Goodman J.B. and Majda A.,
\textit{The validity of the modified equation for nonlinear shock waves,}
J. Comput. Phys. 58 (1985), 336--348.

\bibitem{HayesLeFloch} Hayes B.T. and LeFloch P.G.,
\textit{Nonclassical shocks and kinetic relations~: finite difference schemes,}
SIAM J. Numer. Anal. 35 (1998), 2169--2194.

\bibitem{HouLeFloch} Hou T.Y. and LeFloch P.G.,
\textit{Why nonconservative schemes converge to wrong solutions:
Error analysis,}
Math. of Comput. 62 (1994), 497--530.

\bibitem{JinSlemrod} Jin S. and Slemrod M.,
\textit{Regularization of the Burnett equations via relaxation,}
J. Stat. Phys. 103 (2001), 1009--1033.

\bibitem{KlingenbergLuZhao} Klingenberg, C., Lu, Y. and Zhao, H.,
\textit{$L^1$ singular limit for relaxation and viscosity
approximation of extended traffic models,}
Elect. J. Diff. Equa. 23 (2003), 1--11.

\bibitem{LeFloch} LeFloch P.G.,
{\bf Hyperbolic systems of conservation laws:
The theory of classical and nonclassical shock waves,}
Lectures in Mathematics, ETH Z\"urich, Birkh\"auser, 2002.

\bibitem{Liu} Liu T.P.,
\textit{Hyperbolic conservation laws with relaxation,}
Comm. Math. Phys. 108 (1987), 153--175.

\bibitem{MN} Mascia C. and Natalini R.,
\textit{$L^1$ nonlinear stability of traveling waves for a
hyperbolic system with relaxation,} J. Differential Equations 
132 (1996), 275--292.

\bibitem{Natalini} Natalini R.,
\textit{Recent results on hyperbolic relaxation problems,}
in ``Analysis of systems of conservation laws'' (Aachen, 1997),
Chapman \& Hall/CRC Monogr. Surv. Pure Appl. Math., 99,
Boca Raton, FL, 1999, pp.~128--198.

\bibitem{Slemrod} Slemrod M.,
\textit{Constitutive relations for monatomic gases based on a generalized
rational approximation to the sum of the Chapman-Enskog expansion,}
Arch. Rational Mech. Anal. 150 (1999), 1--22.

\bibitem{Szepessy} Szepessy A.,
\textit{On the stability of Broadwell shocks,} ``Nonlinear
Evolution Partial Differential Equations'', (Beijing 1993),
AMS/IP Stud. Adv. Math. 3, Amer. Math. Soc.,
Providence, RI, 1997, pp. 403--412. 
  
     

\end{thebibliography}
\end{document}